\documentclass{amsart}
\usepackage{amssymb}

\usepackage{graphicx}
\usepackage{amscd}

\newtheorem{theorem}{Theorem}
\theoremstyle{plain}

\newtheorem{corollary}{Corollary}

\newtheorem{lemma}{Lemma}

\newtheorem{proposition}{Proposition}

\numberwithin{equation}{section}
\numberwithin{theorem}{section}
\numberwithin{lemma}{section}
\numberwithin{proposition}{section}
\numberwithin{corollary}{section}

\input{tcilatex}

\begin{document}
\title[Exponential Sums ]{Exponential Sums Along $p-$adic Curves}
\author{W.A. Zuniga-Galindo}
\address{Laboratorio de Computo Especializado, Universidad Autonoma de Bucaramanga,
A.A. 1642, Bucaramanga, Colombia}
\curraddr{Department of Mathematics and Computer Science, Barry University, 11300 N.E.
Second Avenue, Miami Shores, Florida 33161, USA}
\email{wzuniga@mail.barry.edu}
\thanks{Supported by COLCIENCIAS, Grant \# 089-2000}

\begin{abstract}
Let $K$\ be a $p-$adic field, $R$\ the valuation ring of $K$, and \ $P$\ the
maximal ideal of\ $R$. Let $Y\subseteq R^{2}$ be a non-singular closed
curve, and $Y_{m}$ its image in $R/P^{m}\times R/P^{m},$ i.e. \ the reduction%
\textit{\ }modulo\textit{\ }$P^{m}$ of \ $Y$. \ We denote by $\Psi $ an
standard additive character on $K$. In this paper we discuss the estimation
of exponential sums of type$\ \ S_{m}(z,\Psi ,Y,g):=\sum\limits_{x\in
Y_{m}}\Psi (zg(x)),$ with $z\in K$, and $g$\ \ a polynomial function on $Y$.
We show \ that if the p-adic absolute value of $z$ is big enough then the
complex absolute value of $\ \ \ \ S_{m}(z,\Psi ,Y,g)$ is $O(q^{m(1-\beta
(f,g))})$, for a positive constant $\beta (f,g)$ satisfying $0<\beta (f,g)<1$%
.
\end{abstract}

\keywords{exponential sums, $p-$adic fields, algebraic curves}
\subjclass{Primary 11L40, 11T55}
\maketitle

\section{Introduction}

In this paper we shall discuss the estimation of \ exponential sums \ along $%
p-$adic curves. More precisely, let $K$\ be a $p-$adic field, i.e. a finite
algebraic extension of $\mathbb{Q}_{p}$, the field of $p-$adic numbers. Let $%
R$\ be the valuation ring of $K$, \ $P$\ the maximal ideal of\ $R$, and $%
\overline{K}=R/P$\ the residue field \ of $K$. The cardinality of $\overline{%
K}$ is denoted by $q$, thus $\overline{K}=\mathbb{F}_{q}$. For $z\in K$, $%
v(z)\in \mathbb{Z\cup \{+\infty \}}$ denotes \ the valuation of $z$, and $\
\mid z\mid =q^{-v(z)}$ its absolute value. We fix a uniformizing parameter $%
\pi $ for $R$.

We set $X:=R^{n}$ \ as \ $p-$adic analytic variety of dimension $n$, and 
\begin{equation*}
X_{m}:=R/P^{m}\times R/P^{m}\times ..\times R/P^{m},\text{ (}n\text{
factors).}
\end{equation*}

Let $Z\subseteq X$ be a closed subset. We denote by \ $Z_{m}$\ the image of $%
Z$\ in \ $X_{m}$; $Z_{m}$ is named \textit{\ the reduction modulo }$P^{m}$
of \ $Z$.

Let $\Psi $ be an standard additive character on $K$. Thus \ for $z\in K$,

\begin{equation*}
\Psi (z)=\exp (2\pi iT_{r\text{ }K/\mathbb{Q}_{p}}(z)),
\end{equation*}
where $T_{r\text{ }}$ \ denotes \ the trace function.

Let \ \ $Y$ \ be a closed analytic \ subset of $X$, i.e. a closed subset \
which is locally the common locus of zeroes of a finite number of $K$%
-analytic functions. We define 
\begin{equation}
S_{m}(z,\Psi ,Y,g):=\sum\limits_{x\in Y_{m}}\Psi (zg(x)),  \label{eq1}
\end{equation}
where $\ z\in K$, $m=-v(z)\geqq 1$, and $g(x)$ is a non-zero analytic
function on $Y$.

In this paper we discuss the estimation of exponential sums of type (\ref
{eq1}), in the case of $Y$ is a plane non-singular curve. We note  that the
reduction modulo $P$ of $Y$, i.e. $Y_{1}$, may be a singular curve. In this
case, for $\mid z\mid $ big enough, we obtain an \ estimation \ \ of type 
\begin{equation}
\mid S_{m}(z,\Psi ,Y,g)\mid _{\mathbb{C}}\leq A(K,f,g)q^{m(1-\beta (f,g))},
\label{est1}
\end{equation}
where $\mid x\mid _{\mathbb{C}}$ denotes the absolute value of a\ complex
number $x$, and $A(K,f,g)$, $\beta (f,g)$ are\ positive constants, with $%
0<\beta (f,g)<1$ (cf. theorem \ref{theorem1}).

\bigskip

Character sum estimates of Weil-type \ has a long history especially summing
\ over points on a curve defined over a finite field \cite{W}, \cite{B}. A
generalization of the results of the Weil and Bombieri to $p-$adic lifting \
of points over finite fields has been obtained by
Kumar-Helleseth-Calderbrank \cite{KHC}, Li \ \cite{L}, and Voloch-Walker 
\cite{VW1}, \cite{VW2}. More precisely, the \ mentioned works contain
estimates of the form

\begin{equation}
\left| S(z,\Psi ,Y,g)\right| _{\mathbb{C}}\leqq Bp^{m}q^{\frac{1}{2}},
\label{eq2}
\end{equation}
where 
\begin{equation}
S(z,\Psi ,Y,g):=\sum\limits_{x\in Y_{1}}\Psi (zg(x)),
\end{equation}
\ $m=-v(z)$, $q=p^{n}$, and $\ B$ is a constant depending on the function
field of $Y$. The above estimate is good when \ $n$ is large comparing \ $m$%
. \ Our main result considers \ exponential sums along the points of the 
reduction modulo $P^{m}$ of a non-singular curve $Y$ in the case in which \ $%
n$ is fixed and   $m$ is big enough.

In the case of $Y=R^{n}$, Igusa developed a complete theory \ for the
exponential sums \ of type \ (\ref{eq1}) (see e.g. \cite{I}, \cite{D}). \
For other analytic varieties no known results are available, as far as the
author knows.

\section{Preliminary results}

In this section, to seek completeness, we collect some results that it will
be used in the next sections.

\subsection{Serre's measure}

Suppose that $Y\subseteq X$ is non-singular closed analytic subset of
dimension $d$ everywhere. The canonical \ immersion of $Y$ in $X$ induces \
on $Y$ a canonical measure $\mu _{Y},$ completely analog to volume measure
of the real case. If 
\begin{equation*}
I=\{i_{1},i_{2},..,i_{d}\},
\end{equation*}
with \ $i_{1}<i_{2}<..<i_{d}$ is a subset of \ $d$ elements of $\
\{1,2,..,n\}$, we denote by $\omega _{Y,I\text{ }}$ the differential form
induced on $Y$\ by $dx_{i_{1}}\wedge dx_{i_{2}}\wedge ..\wedge dx_{i_{d}}$
and by \ $\alpha _{Y,I\text{ }}$\ the measure corresponding to \ $\omega
_{Y,I\text{ .}}$ Then the canonical measure $\alpha _{Y}$\ is defined as
follows: 
\begin{equation*}
\alpha _{Y}:=\sup_{I}\{\alpha _{Y,I\text{ }}\},
\end{equation*}
where $I$ \ runs through all subsets of $d$\ elements \ of \ $\{1,2,..,n\}$
(cf. \cite{S}).

If \ $Y$ is singular, i.e. if $Y$ is a closed analytic subset of $X$ \ of
dimension \ $\leq d$, we denote by $Y^{reg}$\ the set of points on which $Y$%
\ is smooth of dimension $d$. The canonical measure $\alpha _{Y}$ is defined
on $Y^{reg}$\ \ as a bounded measure (cf. \cite{S}). This allows us \ define
the canonical measure in the singular case by

\begin{equation*}
\alpha _{Y}(A):=\alpha _{Y^{reg}}(A\cap Y^{reg}),\text{ for any }A\subseteq
Y.
\end{equation*}

In \cite{S} Serre \ showed that if $Y$ is closed non-singular analytic
submanifold of dimension $d$ everywhere, then for \ $m$ big enough, it holds
that 
\begin{equation}
\;Card(Y_{m})=\alpha _{Y}(Y)q^{md}.  \label{F}
\end{equation}

Let $x_{0}$ be a fixed point of $Y_{m}$. The proof of (\ref{F}) was
accomplished by showing that 
\begin{equation}
\alpha _{Y}(\{u\in Y\mid u\equiv x_{0}\text{ mod }\pi ^{m}\})=\frac{1}{q^{md}%
},  \label{F1}
\end{equation}
for $m$ big enough (cf. \cite{S}, page 347).

The following proposition follows directly from (\ref{F1}).

\bigskip

\begin{proposition}
\label{prop.1}\ Let $Y\subseteq X$ be a non-singular closed analytic curve, $%
z=\pi ^{-m}u\in K$, $u\in R^{\times }$. Then for $\mid z\mid $big enough, it
holds that
\end{proposition}

\begin{equation}
S_{m}(z,\Psi ,Y,g)=q^{m}\int\limits_{Y}\Psi (zg(x))\alpha _{Y}(x).
\label{R1}
\end{equation}

\bigskip For $Y$ singular, we do not know if there exists a relation \
between $S_{m}(z,\Psi ,Y,g)$ and $\int\limits_{Y}\Psi (zg(x))\alpha _{Y}(x).$

\subsection{Non-archimedean implicit function theorem}

A series 
\begin{equation*}
g(x)=\sum\limits_{i}c_{i}x^{i}\in K[[x]],x=(x_{1},x_{2},..,x_{n})
\end{equation*}
is named \ a \textit{special restricted power series}, abbreviated as SRP,
if $f(0)=0$, and $c_{i}\in P^{\mid i\mid -1}$, $\mid i\mid
=i_{1}+i_{2}+..+i_{n}$, for all $i\neq 0$ in \ $\mathbb{N}^{n}$. This
clearly implies that $f(x)$ in $R[[x]]$. Furthermore $f(x)$ \ is convergent
at every $a\in R^{n}$.

\bigskip

\begin{lemma}[{\protect\cite[Theorem 2.2.1]{I}}]
\label{ift}(Implicit Function Theorem) (i) If 
\begin{equation*}
F_{i}(x,y)\in \ R[[x_{1},..,x_{n},y_{1},..,y_{m}]],
\end{equation*}
\ \ \ $F_{i}(0,0)=0$ for all $i=1,..m,$ $F(x,y)=(F_{1}(x,y),..,F_{m}(x,y))$
\ and further 
\begin{equation*}
\frac{\partial (F_{1},..,F_{m})}{\partial (y_{1},..,y_{m})}(0,0)\not\equiv
{}0\text{ mod }\pi ,
\end{equation*}
in which $\frac{\partial (F_{1},..,F_{m})}{\partial (y_{1},..,y_{m})}$\ is
the Jacobian determinant of the square matrix of size $m\times m$ with $%
\frac{\partial F_{i}}{\partial y_{j}}$\ as its $(i,j)$-entry. Then there
exists a unique $f(x)=(f_{1}(x),..,f_{m}(x))$ with every $f_{i}(x)$ in $%
R[[x]]$ satisfying $f_{i}(0)=0$, and $F(x,f(x))=0$. (ii) If every $F_{i}(x,y)
$\ is an SRP\ \ in $x_{1},..,x_{n},y_{1},..,y_{m}$, then every $f_{i}(x)$ \
is an SRP in $x_{1},..,x_{n}$. Furthermore if $a$ is in $R^{n}$, then $f(a)$
is in $R^{m}$ and $F(a,$ $f(a))=0$, and if $(a,b)$ in $R^{n}\times R^{m}$
satisfies \ $F(a,b)=0$, then $b=f(a)$.
\end{lemma}

\subsection{Exponential sums along parametric curves}

If $A$ \ is open and compact set containing the origin \ and $h:A\rightarrow
R^{2}$ \ is an analytic mapping of the form 
\begin{equation}
h(t)=(t,\eta (t))\text{, with }\eta
(t)=a_{m}t^{m}+\sum\limits_{j=m+1}^{\infty }a_{j}t^{j}\text{, }m\geqq 1\text{%
, and }a_{m}\neq 0\text{.}  \label{PC}
\end{equation}
We call the set 
\begin{equation}
V_{h,A}(R):=\{(t,\eta (t))\mid t\in A\},  \label{PC2}
\end{equation}
a \textit{\ parametric curve} \textit{\ at the origin.}

\bigskip

Given a differential form $dx$ , we denote by $\mid dx\mid $\ the
corresponding Haar measure, normalized such a way that the volume of $R$\ is
1.

If $V_{h,A}(R)$ is a\ parametric curve\ and $g$ a non-zero analytic function
defined on an open set $U$ containing $A$, we define

\begin{equation*}
\mu (h,g):=ord_{t}(g(t,\eta (t))-g\left( 0,0\right) ).
\end{equation*}
If \ $g\left( 0,0\right) =0$, the number $\mu (h,g)$ is the intersection \
multiplicity at \ the origin\ of the analytic curves \ $V_{h,A}(R)$, and $%
V_{g,U}(R)$. In general $\mu (h,g)\geqq mult(h,g)$, where $mult(h,g)$\
denotes the multiplicity of intersection \ of the analytic curves $%
V_{h,A}(R) $, and $V_{g,U}(R)$.

Given a real number $x$, we denote by $[x]$ \ the largest integer satisfying 
$\ [x]\leq $ $x$.

\begin{lemma}
\label{lemma1}Let $V_{h,P^{l}}(R)$ be a parametric curve of type (\ref{PC2}%
), \ $g(x,y)$ a non-constant \ analytic function such that $g(t,\eta
(t))=g(0,0)+t^{\mu (h,g)}\alpha (t)$, $\alpha (0)=c_{0}\in K^{\times }$, and
\ $z=u\pi ^{-m}$\ $\in K,u\in R^{\times }$. \ If $\ l\geqq $ $[\frac{%
m-v(c_{0})}{\mu (h,g)}]+1$, then 
\begin{equation*}
\mid S_{m}(z,\Psi ,V_{h,P^{l}}(R),g)\mid _{\mathbb{C}}\leqq q^{m(1-\frac{1}{%
\mu (h,g)})+\frac{v(c_{0})}{\mu (h,g)}},
\end{equation*}
for $\mid z\mid $\ big enough.
\end{lemma}

\begin{proof}
By proposition \ref{prop.1}, for $\mid z\mid $\ big enough 
\begin{equation*}
S_{m}(z,\Psi ,V_{h,P^{l}}(R),g)
\end{equation*}
can be expressed as an integral \ with respect to Serre's measure as follows:

\begin{equation}
S_{m}(z,\Psi ,V_{h,P^{l}}(R),g)=q^{m}\int\limits_{V_{h,P^{l}}(R)}\Psi
(zg(x,y))\alpha _{V_{h,P^{l}}}.  \label{eq15}
\end{equation}

By using the fact that $V_{h,A}(R)$ is a parametric curve, we \ have that 
\begin{equation}
S_{m}(z,\Psi ,V_{h,P^{l}}(R),g)=q^{m}\int\limits_{P^{l}}\Psi (zg(t,\eta
(t)))\mid dt\mid .  \label{eq17}
\end{equation}

\ 

We put \ 
\begin{equation}
g(t,\eta (t))=g(0,0)+t^{\mu (h,g)}\alpha (t),  \label{eq19}
\end{equation}
with $\alpha (0)=c_{0}\in K^{\times }$.

The integral in (\ref{eq17}) admits the following expansion:

\begin{equation}
S_{m}(z,\Psi ,V_{h,P^{l}}(R),g)=q^{m}\Psi (zg(0,0))\sum\limits_{j=l}^{\infty
}q^{-j}\int\limits_{R^{\times }}\Psi (z\pi ^{j\mu (h,g)}\alpha (\pi
^{j}t))\mid dt\mid .  \label{eq21}
\end{equation}

If \ $l\geqq \max $ $[\frac{m-v(c_{0})}{\mu (h,g)}]+1\geqq v(c_{0})+1$, it
holds that 
\begin{equation}
v(\alpha (\pi ^{j}t))=v\left( c_{0}\right) \text{, for every \ }t\in
R^{\times }\text{, and \ \ }j\geqq l\text{,}  \label{eq22a}
\end{equation}

and 
\begin{equation}
\Psi (z\pi ^{j\mu (h,g)}\alpha (\pi ^{j}t))=1\text{, for every \ }t\in
R^{\times }\text{, and \ \ }j\geqq l\text{.}  \label{eq22}
\end{equation}

Therefore from (\ref{eq21}), (\ref{eq22a}), and (\ref{eq22}), it follows that

\begin{equation}
\mid S_{m}(z,\Psi ,V_{h,P^{l}}(R),g)\mid _{\mathbb{C}}=q^{m-l}\leqq q^{m(1-%
\frac{1}{\mu (h,g)})+\frac{v(c_{0})}{\mu (h,g)}}.  \label{eq23}
\end{equation}

\qquad \qquad \qquad 
\end{proof}

\bigskip

\section{Exponential sums along non-singular curves}

In this section we \ discuss the estimation of exponential sums of type (\ref
{eq1})\ along non-singular curves.

\bigskip Let $f:U\rightarrow K$ be \ an analytic function on an open and
compact neighborhood $U$ $\subseteq K^{2}$ of a point $(x_{0},y_{0})\in
K^{2} $. By \ an $K-$analytic curve at $(x_{0},y_{0})$, we mean an analytic\
set of the form 
\begin{equation}
V_{f,U}(K)=\{(x,y)\in U\mid f(x,y)=0\}.  \label{eq3}
\end{equation}

We set $V_{f,U}(R):=$\ $V_{f,U}(K)\cap R^{2}$.

\bigskip

Let $W$\ be an open and compact \ set containing the origin, and $V_{f,W}(K)$
\ an analytic curve at the origin, such that the origin is a smooth point,
i.e. 
\begin{equation}
f(x,y)=ax+by+(\text{higher order terms}),\text{ \ }a\neq 0\text{ or \ }b\neq
0.  \label{eq27}
\end{equation}
Since $f(x,y)$\ is a convergent series, by multiplying it \ by a non-zero
constant $c_{0}\in R$, we may suppose that $f(x,y)\in R[[x,y]]\setminus
P[[x,y]]$.

\bigskip

We define 
\begin{equation}
L(f,(0,0)):=\min \{v(a),v(b)\},  \label{eq28}
\end{equation}
with $a$, $b$ as in (\ref{eq27}). The Jacobian criteria implies that \ $%
L(f,(0,0))=0$ if and only if \ the origin is a non-singular point of the
reduction modulo $\pi $ of $f(x,y)$. If $L(f,(0,0))=v(b)\neq 0$, \ by making
a change of coordinates of the form \ $x=\pi ^{v(b)+1}x^{\prime },$ $y=\pi
^{v(b)+1}y^{\prime }$, we get that \ $f(\pi ^{v(b)+1}x^{\prime },\pi
^{v(b)+1}y^{\prime })=\pi ^{2v(b)+1}\widetilde{f}(x^{\prime },y^{\prime })$
with $L(\widetilde{f},(0,0))=0.$ Furthermore $\widetilde{f}(x^{\prime
},y^{\prime })$ is an SRP.

\bigskip

\begin{lemma}
\label{lemma2}Let $V_{f,P^{l}\times P^{l}}(R)$ be $\ $an analytic curve at
the origin, such that the origin is a smooth point, and $g$ an analytic
function such that $g\mid _{V_{f,P^{l}\times P^{l}}(R)}$ $\neq 0$, with $%
g(0,0)=0.$ There exist constants $C(K,f,g)$, $\gamma (f,g)$ such that if $%
l\geqq \lbrack \frac{m-\gamma (f,g)}{\mu (f,g)}]+1$, then 
\begin{equation*}
\mid S_{m}(z,\Psi ,V_{f,P^{l}\times P^{l}}(R),g)\mid _{\mathbb{C}}\leqq
C(f,g)q^{m(1-\frac{1}{\mu (f,g)})},
\end{equation*}
for $\mid z\mid $\ big enough.
\end{lemma}

\begin{proof}
We assume without loss of generality that $f(x,y)\in R[[x,y]]\setminus
P[[x,y]]$,\ and that $L(f,(0,0))=v(b)\geqq 1$. Suppose that and 
\begin{equation*}
g(x,y)=g_{D}(x,y)+(\text{higher order terms}),\text{ }
\end{equation*}
with $g_{D}(x,y)$ a $\ $homogeneous polynomial of \ degree $D$. By
proposition \ref{prop.1}, for $\mid z\mid $\ big enough 
\begin{equation*}
S_{m}(z,\Psi ,V_{f,P^{l}\times P^{l}}(R),g)
\end{equation*}
can be expressed as an integral \ with respect to Serre's measure as follows:

\begin{equation}
S_{m}(z,\Psi ,V_{f,P^{l}\times
P^{l}}(R),g)=q^{m}\int\limits_{V_{f,P^{l}}(R)}\Psi (zg(x,y))\alpha
_{V_{f,P^{l}}}.  \label{eq29a}
\end{equation}

By making a change of coordinates of the form \ $x=\pi ^{v(b)+1}x^{\prime },$
$y=\pi ^{v(b)+1}y^{\prime }$ in (\ref{eq29a}), and assuming that \ $l\geqq
m+1\geqq L(f,(0,0)+1$, \ we have that 
\begin{equation}
S_{m}(z,\Psi ,V_{f,P^{l}\times P^{l}}(R),g)=q^{-v(b)-1}S_{m}(z\pi ^{D\left(
v(b)+1\right) },V_{f^{\ast },P^{l-v(b)-1}\times P^{l-v(b)-1}}(R),g^{\ast }),
\label{eq29}
\end{equation}
where 
\begin{equation*}
f^{\ast }(x^{\prime },y^{\prime })=\pi ^{-2v(b)+1}f(\pi ^{v(b)+1}x^{\prime
},\pi ^{v(b)+1}y^{\prime }),
\end{equation*}
is an SRP, and 
\begin{equation*}
g^{\ast }(x^{\prime },y^{\prime })=\pi ^{-D\left( v(b)+1\right) }g(\pi
^{v(b)+1}x^{\prime },\pi ^{v(b)+1}y^{\prime }).
\end{equation*}

Since 
\begin{equation*}
L(f^{\ast },(0,0))=0,
\end{equation*}
it holds that 
\begin{equation*}
\frac{\partial f^{\ast }(x^{\prime },y^{\prime })}{\partial y^{\prime }}%
(0,0)\not\equiv {}0\text{ mod }\pi .
\end{equation*}
Thus, \ by implicit function theorem (see lemma \ref{ift}) $V_{f^{\ast
},P^{l-v(b)-1}\times P^{l-v(b)-1}}(R)$ is a parametric curve, i.e. there
exists \ an SRP function $h(t)$ such that 
\begin{equation}
V_{f^{\ast },P^{l-v(b)-1}\times P^{l-v(b)-1}}(R)=\{(t,h(t))\mid t\in
P^{l-v(b)-1}\}.  \label{eq30}
\end{equation}
Now, we set $g^{\ast }\left( h(t)\right) =t^{\mu (f,g)}\beta (t)$, with $%
\beta (0):=\gamma (f,g)\in K^{\times }$. The result follows from (\ref{eq29}%
), and (\ref{eq30}) by lemma \ref{lemma1}.\bigskip 
\end{proof}

\bigskip We note that the previous lemma is valid if $g(0,0)\neq 0$. Indeed,
if $g(x,y)=g(0,0)+g^{\prime }(x,y),$ with $g^{\prime }(0,0)=0$, then 
\begin{equation*}
S_{m}(z,\Psi ,V_{f,P^{l}\times P^{l}}(R),g)=\Psi (zg(0,0))S_{m}(z,\Psi
,V_{f,P^{l}\times P^{l}}(R),g^{\prime }).
\end{equation*}

\bigskip

Given a non-constant polynomial $f(x,y)\in R[x,y]$, we define \ for every
non-singular point $P\in V_{f}(R)=\{(x,y)\in R^{2}\mid f(x,y)=0\},$ the
number 
\begin{equation*}
L(f,P):=\min \{v(\frac{\partial f}{\partial x}(P)),v(\frac{\partial f}{%
\partial y}(P))\}.
\end{equation*}

This definition \ generalizes that given in (\ref{eq28}).

\begin{proposition}
\label{prop.3}Let $f(x,y)\in R[x,y]$ be a non-constant polynomial, \ such
that $V_{f}(R)$ does not have singular points on $R^{2}$. There is a
positive constant $c(f)$, depending only on $f$ \ such that 
\begin{equation*}
L(f,P)\leqq c(f),\,\,\,\text{for all}\,\,\ P\in V_{f}(R).
\end{equation*}
\end{proposition}

\begin{proof}
\bigskip\ By contradiction, we suppose that $L(f,P)$ is not bounded on $%
V_{f}(R)$. So there is a sequence $(Q_{i})_{i\in \mathbb{N}}$ of points of $%
\ V_{f}(R)$ satisfying lim $L(f,Q_{i})\,\,\longrightarrow \infty $, when $%
i\longrightarrow \infty $. Since $\ V_{f}(R)$ is a compact set, the sequence 
$Q_{i}\,\,$ has a limit point $\overline{Q}\in V_{f}(R)\,$. But $\overline{Q}
$ is \ a singular point of $V_{f}(R)$, contradiction. \quad 
\end{proof}

\bigskip The numbers $L(f,P)$\ were introduced by A. N\'{e}ron and they
appear naturally in the computation of several \ p-adic integrals (see \cite
{Z-G1}, \cite{Z-G2} and the references therein).

\bigskip

Given two non-zero polynomial functions $f(x,y),g(x,y)$ with coefficients in 
$R$, such that $g(0,0)=f(0,0)=0$, we define 
\begin{equation}
\sigma _{f}(g):=\sup_{P\in V_{f}(R)}\{mult_{P}(f,g)\}\in \mathbb{N}\setminus
\left\{ 0\right\} .  \label{s}
\end{equation}

\begin{theorem}
\label{theorem1}Let $f(x,y),g(x,y)\in R[x,y]$, be non-zero polynomials, with 
$g(0,0)=f(0,0)=0$. If the curve $V_{f}(R)$ is non-singular, then \ there
exists a constant $A(K,f,g)$ such that 
\begin{equation*}
\mid S_{m}(z,\Psi ,V_{f}(R),g)\mid _{\mathbb{C}}\leqq A(K,f,g)q^{m(1-\frac{1%
}{\sigma _{f}(g)})},
\end{equation*}
for $\mid z\mid $ big enough.
\end{theorem}

\begin{proof}
\bigskip We set a positive integer $l=m+1\geqq c(f)+1$, where $c(f)$ is the
constant whose existence was established in proposition \ref{prop.3}.

The compact set $V_{f}(R)$ admits the following covering: 
\begin{equation}
V_{f}(R)=\bigcup\limits_{i=0}^{r}V_{f,D_{i}}(R),  \label{cover}
\end{equation}
with 
\begin{equation}
V_{f,D_{i}}(R):=V_{f}(R)\cap (x_{i},y_{i})+P^{l}\times P^{l},
\end{equation}
and $Q_{i}=(x_{i},y_{i})\in V_{f}(R),$ $i=1,..,r$. Because each $V_{f,D_{i}}$
$(R)$ \ is an analytic curve \ at $(x_{i},y_{i})$, \ it follows from
covering (\ref{cover}) that

\begin{equation}
\mid S_{m}(z,\Psi ,V_{f}(R),g)\mid _{\mathbb{C}}\leqq
\sum\limits_{i=1}^{r}\mid S_{m}(z,\Psi ,V_{f},_{D_{i}}(R),g)\mid _{\mathbb{C}%
}.  \label{eq32}
\end{equation}

By \ applying lemma \ref{lemma2} \ in (\ref{eq32}), we have 
\begin{equation}
\mid S_{m}(z,\Psi ,V_{f}(R),g)\mid _{\mathbb{C}}\leqq
\sum\limits_{i=1}^{r}C_{i}(f,g,K)q^{m(1-\frac{1}{\mu _{Q_{i}}\left(
f,g\right) })}.  \label{eq34}
\end{equation}

The result follows from (\ref{eq34}) by observing that $\max_{i}\left\{ 1-%
\frac{1}{\mu _{Q_{i}}\left( f,g\right) }\right\} \leqq (1-\frac{1}{\sigma
_{f}(g)})$.
\end{proof}

Let 
\begin{equation*}
f(x)=\alpha _{0}\left( \prod\limits_{i=1}^{r}(x-\alpha _{i})^{e_{i}}\right)
Q(x)\in K[x]
\end{equation*}
be a non-constant polynomial in one variable, \ and \ $Q(x)$ an irreducible
polynomial of degree greater or equal two.

We put

\begin{equation}
\sigma (f):=\max_{i}\{e_{i}\}.
\end{equation}

The following corollary follows immediately from theorem \ref{theorem1}.

\begin{corollary}
Let $f(x)$ be a non-constant polynomial in one variable. Then with the above
notation, it holds for $\mid z\mid $ big enough that 
\begin{equation}
\mid \sum\limits_{x\in R_{m}}\Psi (zf(x))\mid _{\mathbb{C}}\leqq
A(K,f,g)q^{m(1-\frac{1}{\sigma (f)})}.  \label{e2}
\end{equation}
\end{corollary}


\bigskip

\begin{thebibliography}{99}
\bibitem{B}  Bombieri E., On exponential sums in finite fields, Amer. J.
Math. 88 (1966), 71-105.

\bibitem{D}  Denef J., Report on Igusa's local zeta functions, Seminaire
Bourbaki 741 1990-1991.

\bibitem{I}  Igusa J. -I., An introduction \ to the theory of local zeta
functions, AMS/IP Studies in Advanced Mathematics, vol. 14, 2000.

\bibitem{KHC}  Kumar P.V., Helleseth T., and Calderbrank A. R., \ An upper
bound for some exponential sums over Galois rings and applications, IEEE
Trans. Info. Theory 41 (1995) 456-468.

\bibitem{L}  Li Wen-Ching Winnie, Character sums \ over $p-$adic fields, J.
Number Theory 74, No.2, 181-229 (1999).

\bibitem{S}  Serre J-P., Quelques applications du th\'{e}or\`{e}me de densit%
\'{e} de Cebotarev, Pub. I.H.E.S, vol. 54,1981.

\bibitem{VW1}  Voloch J.F, Walker J.L, Euclidean weights of codes from
elliptic curves over rings, Trans. Amer. Math. Soc. 352, No.11, 5063-5076
(2000).

\bibitem{VW2}  Voloch J.F, Walker J.L, Codes over rings from curves of
higher genus, IEEE Trans. Info. Theory. 45, No. 6, 1768-1776 (1999).

\bibitem{W}  Weil A.,  On some exponential sums, Proc. National Academy of
Science 34 (1948), 204-207.

\bibitem{Z-G1}  Zuniga-Galindo W.A., Igusa's local zeta functions of
semiquasihomogeneous polynomials, Trans. Amer. Math. Soc.353, (2001),
3193-3207.

\bibitem{Z-G2}  Zuniga-Galindo W.A., Local zeta functions and Newton
polyhedra, preprint 2000.
\end{thebibliography}
\end{document}